\documentclass{amsart}
\usepackage{amssymb,amsthm,comment,color,enumerate,adjustbox,enumitem,fancyhdr}
\usepackage[outer=1.5in, inner=1.25in, top=1.25in, bottom=1.25in]{geometry}
\usepackage[english]{babel}

\newtheorem{teo}{Theorem}[]
\newtheorem{cor}[teo]{Corollary}

\begin{document}
	
\title{On the Cayleyness of Praeger-Xu graphs}

\author{Marco Barbieri}
\address{Dipartimento di Matematica ``Felice Casorati", University of Pavia, Via Ferrata 5, 27100 Pavia, Italy} 
\email{marco.barbieri07@universitadipavia.it}
\author{Valentina Grazian}
\address{Dipartimento di Matematica e Applicazioni, University of Milano-Bicocca, Via Cozzi 55, 20125 Milano, Italy} 
\email{valentina.grazian@unimib.it}
\author{Pablo Spiga}
\address{Dipartimento di Matematica e Applicazioni, University of Milano-Bicocca, Via Cozzi 55, 20125 Milano, Italy} 
\email{pablo.spiga@unimib.it}

\begin{abstract}
We give a sufficient and necessary condition for a Praeger-Xu graph to be a Cayley graph. 
\end{abstract}


\maketitle	
\section{Scope of this note}
The Praeger-Xu graphs, introduced by Praeger and Xu in \cite{PX.TwicePrimeSymmetric}, have exponentially large groups of automorphisms, with respect to the number of vertices. This fact causes various complications with regard to many natural questions.

In their recent work \cite{JPS.CayleynessPrXu}, Jajcay, Poto\v{c}nik and Wilson gave a sufficient and necessary condition for a Praeger-Xu graph to be a Cayley graph. Explicitly, \textup{\cite[Theorem 1.1]{JPS.CayleynessPrXu}} states that, for any positive integer $n\geq 3$, $n\ne 4$, and for any positive integer $k\leq n-1$, the Praeger-Xu graph $\textup{PX}(n,k)$ is a Cayley graph if, and only if, one of the following holds
\begin{enumerate}
	\item[(i)] the polynomial $t^n+1$ has a divisor of degree $n-k$ in $\mathbb{Z}_2[t]$;
	\item[(ii)] $n$ is even, and there exist polynomials $f_1,f_2,g_1,g_2,u,v\in \mathbb{Z}_2[t]$ such that $u,v$ are palindromic of degree $n-k$, and
	\begin{equation}\label{eq:JPS}
		t^n+1 = f_1(t^2)u(t) + tg_1(t^2)v(t) = f_2(t^2)v(t) + tg_2(t^2)u(t).
	\end{equation}
\end{enumerate}

\par Our aim here is to prove that  (ii) implies (i), thus obtaining the following refinement. (It can be verified that $\textup{PX}(4,1)$, $\textup{PX}(4,2)$ and $\textup{PX}(4,3)$ are Cayley graphs.)
\begin{teo}\label{th2}
	For any positive integer $n\geq 3$ and for any positive integer $k\leq n-1$, the Praeger-Xu graph $\textup{PX}(n,k)$ is a Cayley graph if, and only if, the polynomial $t^n+1$ has a divisor of degree $n-k$ in $\mathbb{Z}_2[t]$.
\end{teo}
\par Using the factorization of  $t^n+1$ in $\mathbb{Z}_2[t]$, we give a purely  arithmetic condition for the Cayleyness of $\textup{PX}(n,k)$.
Let $\varphi$ be the Euler $\varphi$-function and, for every positive integer $d$, let
\begin{equation*}
	\omega(d):=\min \left\lbrace c\in \mathbb{N} \mid d \text{ divides } 2^c-1 \right\rbrace
\end{equation*}
be the multiplicative order of $2$ modulo $d$.
\begin{cor}\label{th3}
Let $a$ be a non-negative integer, let $b$ be an odd positive integer, let $n:=2^ab$ with $n\geq 3$, and let $k$ be a positive integer with $k\leq n-1$. The Praeger-Xu graph $\textup{PX}(n,k)$ is a Cayley graph if, and only if, $k$ can be written as
	\begin{equation}\label{eq:k}
		k= \sum_{d\mid b} \alpha_d \omega(d), \;\;\;\text{for some integers \;}\alpha_d\text{ with \;} 0\leq \alpha_d \leq \frac{2^a\varphi(d)}{\omega(d)}.
	\end{equation}
\end{cor}

\section{Proof of Theorem \ref{th2}}
Suppose (ii) holds. We aim to show that $t^n+1$ is divisible by a polynomial of degree $n-k$ in $\mathbb{Z}_2[t]$, implying (i). Working in characteristic $2$,~\eqref{eq:JPS} can be written as
\begin{equation*}
	t^n+1=f_1^2(t)u(t)+tg_1^2(t)v(t)=f_2^2(t)v(t)+tg_2^2(t)u(t),
\end{equation*}
in short
\begin{equation}\label{eq:1}
	t^n+1=f_1^2u+tg_1^2v=f_2^2v+tg_2^2u.
\end{equation}

If $g_1=0$ or if $g_2=0$, then the result follows from~\eqref{eq:1} and the fact that $u$ and $v$ have degree $n-k$. Therefore, for the rest of the argument we may suppose that $g_1,g_2\ne 0$. Moreover, observe that $f_1,f_2\ne 0$, because $t$ does not divide $t^n+1$.
\par We introduce four polynomials $u_e,u_o,v_e,v_o \in \mathbb{Z}_2[t]$ such that
	\[ u:=u_e^2+tu_o^2 , \qquad \qquad	v:=v_e^2+tv_o^2.\]
Substituting these expansions for $u$ and $v$ in~\eqref{eq:1}, we get
\begin{align*}
	t^n+1&=f_1^2u_e^2+t^2g_1^2v_o^2+t(f_1^2u_o^2+g_1^2v_e^2),\\
	t^n+1&=f_2^2v_e^2+t^2g_2^2u_o^2+t(f_2^2v_o^2+g_2^2u_e^2).
\end{align*}
Recall that $n$ is even. By splitting the equalities in even and odd degree terms, we obtain
\begin{align*}
	t^n+1&=f_1^2u_e^2+t^2g_1^2v_o^2, \qquad \qquad 0=t(f_1^2u_o^2+g_1^2v_e^2),\\
	t^n+1&=f_2^2v_e^2+t^2g_2^2u_o^2, \qquad \qquad 	0=t(f_2^2v_o^2+g_2^2u_e^2).
\end{align*}
Set $m:=n/2$. Since we are working in characteristic $2$, we get
\begin{alignat}{2}\label{eq:21}
	&t^m+1=f_1u_e+tg_1v_o, \qquad \qquad &&t^m+1=f_2v_e+tg_2u_o,\\
\label{eq:23}
	&f_1u_o=g_1v_e,  &&f_2v_o=g_2u_e.
\end{alignat}

Since $u$ and $v$ are palindromic by assumption, we get $1=u(0) =u_e(0)$ and $1=v(0) =v_e(0)$. In particular both $u_e$ and $v_e$ are not zero. From~\eqref{eq:21} and~\eqref{eq:23}, we obtain
\begin{equation}\label{eq:3}
\begin{split}
	f_1=&\frac{t^m+1}{u_ev_e+tu_ov_o}v_e, \qquad \qquad g_1=\frac{t^m+1}{u_ev_e+tu_ov_o}u_o,\\
	f_2=&\frac{t^m+1}{u_ev_e+tu_ov_o}u_e,	\qquad \qquad g_2=\frac{t^m+1}{u_ev_e+tu_ov_o}v_o.
\end{split}
\end{equation}

Our candidate for the searched divisor of  $t^n + 1$ is $s:=u_ev_e+tu_ov_o$. Let us show first that $\deg(s)=n - k$. 
Since $u_ev_e$ and $u_ov_o$ have even degree, we deduce 
\[\deg(s) = \max \{\deg(u_ev_e), \deg(tu_ov_o)\}.\]
Recall $u=u_e^2+tu_o^2$ and $v=v_e^2+tv_o^2$. If $n-k$ is even, then 
\begin{equation*}
	\deg(u_e)=\deg(v_e)=\frac{n - k}{2}, \;\;\textup{and}\;\; \deg(u_o)=\deg(v_o)< \frac{n-k-1}{2}.
\end{equation*}
On the other hand, if $n-k$ is odd, then
\begin{equation*}
	\deg(u_e)=\deg(v_e)<\frac{n-k}{2}, \;\;\textup{and}\;\; \deg(u_o)=\deg(v_o) = \frac{n-k-1}{2}.
\end{equation*}
Therefore, in both cases, $\deg(s) = n-k$.

It remains to prove that $s$ divides $t^n+1$. Since $f_1,g_1,f_2,g_2$ are polynomials, by~\eqref{eq:3},  $s$ divides
\begin{equation*}
	\gcd((t^m+1)v_e,(t^m+1)v_o,(t^m+1)u_e,(t^m+1)u_o)=(t^m+1)\gcd(v_e,v_o,u_e,u_o).
\end{equation*} 
Observe that $\gcd(v_e,v_o,u_e,u_o)$ divides $f_1u_e+tg_1v_o$, hence, in view of the first equation in~\eqref{eq:21},  $\gcd(v_e,v_o,u_e,u_o)$ divides $t^m+1$. Therefore, $s$ divides $(t^m+1)^2=t^n+1$.

\section{Proof of Corollary \ref{th3}}
\label{sec3}
By Theorem \ref{th2}, deciding if a Praeger-Xu graph $\textup{PX}(n,k)$ is a Cayley graph is tantamount deciding if $t^n+1$ admits a divisor of order $k$ in $\mathbb{Z}_2[t]$. An immediate way to proceed is to study how $t^n+1$ can be factorized in irreducible polynomials.

Let $n=2^ab$, with $\gcd(2,b)=1$. Since we are in characteristic $2$ we have
\begin{equation*}
	t^n+1 = t^{2^ab}+1 = \left( t^b+1 \right)^{2^a}.
\end{equation*}
Furthermore, if $\lambda_d(t)\in\mathbb{Z}[t]$ denotes the $d$-th cyclotomic polynomial, then 
\begin{equation*}
	t^b+1 = \prod_{d|b} \lambda_d(t)
\end{equation*} 
is the factorization of $t^b+1$ in irreducible polynomials over $\mathbb{Q}[t]$, by Gauss' theorem.
Since the Galois group of any field extension of $\mathbb{Z}_2$ is a cyclic group generated by the Frobenius automorphism, the degree of an irreducible factor of $\lambda_d(t)$ in $\mathbb{Z}_2[t]$ is the smallest $c$ such that a $d$-th primitive root $\zeta$ elevated to $2^c$ is $\zeta$, that is, $\omega(d)$. Hence $\lambda_d(t)$ in $\mathbb{Z}_2[t]$ factorizes into $\varphi(d)/\omega(d)$ irreducible polynomials, each having  degree $\omega(d)$.

Therefore, $t^n+1\in\mathbb{Z}_2[t]$ has a divisor of degree $k$ if, and only if, $k$ can be written as the sum of some $\omega(d)$'s, each summand repeated at most $2^a\varphi(d)/\omega(d)$ times, which is exactly~\eqref{eq:k}.

\bibliographystyle{alpha}
\bibliography{CPX.bib}

\end{document}